\input amstex
\documentstyle{amsppt}

\def \wt{\widetilde}
\def \ot{\otimes}
\def \into{\hookrightarrow}
\def \si{\sigma}
\def \de{\delta}
\def \ga{\gamma}
\def \ba{\bold a}
\def \bb{\bold b}
\def \fa{\frak a}
\def \fb{\frak b}

\def \Ext{\operatorname{Ext}}
\def \Tor{\operatorname{Tor}}
\def \H{\operatorname{H}}
\def \Hom{\operatorname{Hom}}
\def \op{\operatorname{op}}

\topmatter
\title
        Hochschild Cohomology of Triangular Matrix Algebras
\endtitle

\author
       Jorge A. Guccione and Juan J. Guccione
\endauthor

\address
     Jorge Alberto Guccione, Departamento de Matem\'atica, Facultad de
     Ciencias Exactas y Naturales, Pabell\'on 1 - Ciudad Universitaria,
     (1428) Buenos Aires, Argentina.
\endaddress

\email
     vander\@dm.uba.ar
\endemail

\address
     Juan Jos\'e Guccione, Departamento de Matem\'atica, Facultad de Ciencias
     Exactas y Naturales, Pabell\'on 1 - Ciudad Universitaria, (1428)
     Buenos Aires, Argentina.
\endaddress

\email
     jjgucci\@dm.uba.ar
\endemail

\abstract
Let $E = \Bigl(\smallmatrix A & M \\ 0 & B\endsmallmatrix \Bigr)$ be a
triangular algebra, where $A$ and $B$ are algebras over an arbitrary
commutative ring $k$ and $M$ is an $(A,B)$-bimodule. We prove the existence
of two long exact sequence of $k$-modules relating the Hochschild
cohomology of $A$, $B$ and $E$. We also study the structure of the maps of
the first of these exact sequences.
\endabstract

\subjclass\nofrills{{\rm 2000} {\it Mathematics Subject
Classification}.\usualspace} Primary 16E40; Secondary 16S50 \endsubjclass

\thanks
Supported by UBACYT 01/TW79 and CONICET
\endthanks

\endtopmatter

\document

\head Introduction \endhead

Let $k$ be an arbitrary commutative ring with unit, $A$ and $B$ two
$k$-algebras with unit, $M$ an $(A,B)$-bimodule, $E = \Bigl(\smallmatrix A
& M \\ 0 & B\endsmallmatrix \Bigr)$ the triangular algebra and $X$ an
$E$-bimodule. Let $1_A$ and $1_B$ be the unit elements of $A$ and $B$
respectively. Let us write $X_{AA} = 1_AX1_A$, $X_{AB} = 1_AX1_B$, $X_{BA}
= 1_BX1_A$ and $X_{BB} = 1_BX1_B$. For example, for $X = E$, we have
$X_{AA} = A$, $X_{AB} = M$, $X_{BA} = 0$ and $X_{BB} = B$.

The purpose of this paper is to prove the following results:

\proclaim{Theorem 1} There exists a long exact sequence
$$
\align
0 & @>>> \H^0(E,X) @>j^0>> \H^0(A,X_{AA})\oplus \H^0(B,X_{BB}) @>\de^0>>
\Ext^0_{A\ot B^{\op},k}(M,X_{AB}) @>>> \\
& @>\pi^0>> \H^1(E,X) @>j^1>> \H^1(A,X_{AA})\oplus \H^1(B,X_{BB})
@>\de^1>> \Ext^1_{A\ot B^{\op},k}(M,X_{AB}) @>>>\dots,
\endalign
$$
where $\Ext^*_{A\ot B^{\op},k}(M,X_{AB})$ denotes the $\Ext$ groups of the
$A\ot B^{\op}$-module $M$, relative to the family of the $A\ot
B^{\op}$-linear epimorphisms which split as $k$-linear morphisms.
\endproclaim

Let $\pi\:E\to B$ the ring morphism defined by $\pi\Bigl(\Bigl(
\smallmatrix a & m \\ 0 & b\endsmallmatrix \Bigr)\Bigr) = b$. We let $B_E$
denote the ring $B$ consider as an $E$-bimodule via $\pi$.

\proclaim{Theorem 2} There exists a long exact sequence
$$
\align
0 & @>>> \Ext^0_{E\ot B^{\op},k}(B_E,X1_B) @>>> \H^0(E,X) @>>>
\H^0(A,X_{AA}) @>>>  \\
& @>>> \Ext^1_{E\ot B^{\op},k}(B_E,X1_B) @>>> \H^1(E,X) @>>>
\H^1(A,X_{AA}) @>>> \dots,
\endalign
$$
where $\Ext^*_{E\ot B^{\op},k}(B_E,X1_B)$ denotes the $\Ext$ groups of the
$E\ot B^{\op}$-module $B_E$, relative to the family of the $E\ot
B^{\op}$-linear epimorphisms which split as $k$-linear morphisms.
\endproclaim

Let $C$ be a $k$-algebra. The complex $\Hom_{C^e}((C^{*+2},b'_*),C)$ is a
differential graded algebra via the cup product
$$
(fg)(x_0\ot\cdots\ot x_{r+s+1}) = f(x_0\ot\cdots\ot x_r\ot 1_E)g(1_E\ot
x_{r+1}\ot\cdots\ot x_{r+s+1}),
$$
where $f\in \Hom_{C^e}(C^{r+2},C)$ and $g\in \Hom_{C^e}(C^{s+2},C)$. Hence
$\H^*(E,E)$ becomes a graded associative algebra. As it is well known,
$\H^*(E,E)$ is a graded commutative algebra. Similarly, for a left
$C$-module $Y$, the complex $\Hom_C((C^{*+1}\ot Y,b'_*),Y)$ is a
differential graded algebra via
$$
(fg)(x_0\ot\cdots\ot x_{r+s}\ot y) = f(x_0\ot\cdots \ot x_r\ot g(1\ot
x_{r+1}\ot\cdots \ot x_{r+s}\ot y),
$$
where $f\in \Hom_C(C^{r+1}\ot Y,Y)$ and $g\in \Hom_C(C^{s+1}\ot Y,Y)$. Hence
$\Ext_{C,k}^*(Y,Y)$ is a graded associative algebra. In the following
theorem we consider $\H^*(A,A)$, $\H^*(B,B)$, $\H^*(E,E)$ and $\Ext^*_{A\ot
B^{\op},k}(M,M)$ equipped with these algebra structures.

\proclaim{Theorem 3} Assume that $X = E$. Then

\roster

\item The map $\H^*(E,E) @>j^*>> \H^*(A,A) \oplus \H^*(B,B)$ is a
morphism of graded rings.

\item The maps $\H^*(A,A) \into \H^*(A,A) \oplus \H^*(B,B) @>\de^*>>
\Ext^*_{A\ot B^{\op},k}(M,M)$ and $\H^*(B,B) \into \H^*(A,A) \oplus
\H^*(B,B) @>\de^*>> \Ext^*_{A\ot B^{\op},k}(M,M)$ are morphisms of graded
rings.

\item The map $\Ext^*_{A\ot B^{\op},k}(M,M) @>\pi^*>> \H^{*+1}(E,E)$ has
image in the annihilator of $\bigoplus_{n\ge 1} \H^n(E,E)$.

\endroster

\endproclaim

Theorems~1 and 2, which generalize a previous result of \cite{H}, were
established in \cite{M-P} under the assumptions that $k$ is a field, $A$
and $B$ are finite dimensional $k$-algebras, $M$ is a finitely generated
$(A,B)$-bimodule and $X = E$. As was pointed out in \cite{M-P}, this
version of Theorem~1, also follows from a result of \cite{C}. Theorem~3 was
proved in \cite{G-M-S, Setion 5}, under the assumptions that $k$ is a
field, $A$ is a finite dimensional $k$-algebra and $E$ is a one point
extension of $A$.

Our proofs are elementary. The main tool that we use is the existence of a
simple relative projective resolution of $E$.

Next, we enunciate the homological versions of Theorems~1 and 2. Similar
methods to the ones used to prove Theorems~1 and 2 work in the homological
context. We left the task of giving the proofs to the reader.

\proclaim{Theorem 1'} Let $X_{BA} = 1_BX1_A$. There exists a long exact
sequence
$$
\align
\dots & @>>> \Tor_1^{A\ot B^{\op},k}(M,X_{BA}) @>>> \H_1(A,X_{AA})\oplus
\H_1(B,X_{BB}) @>>> \H_1(E,X) @>>> \\
& @>>> \Tor_0^{A\ot B^{\op},k}(M,X_{BA}) @>>>  \H_0(A,X_{AA})\oplus
\H_0(B,X_{BB}) @>>> \H_0(E,X) @>>> 0,
\endalign
$$
where $\Tor_*^{A\ot B^{\op},k}(M,X_{BA})$ denotes the $\Tor$ groups of the
$A\ot B^{\op}$-module $M$, relative to the family of the $A\ot
B^{\op}$-linear epimorphisms which split as $k$-linear morphisms.
\endproclaim

\proclaim{Theorem 2'} There exists a long exact sequence
$$
\align
\dots & @>>> \H_1(A,X_{AA}) @>>> \H_1(E,X) @>>> \Tor_1^{E\ot
B^{\op},k}(B_E,1_BX) @>>>\\
& @>>> \H_0(A,X_{AA}) @>>> \H_0(E,X) @>>> \Tor_0^{E\ot B^{\op},k}(B_E,1_BX)
@>>> 0,
\endalign
$$
where $\Tor_*^{E\ot B^{\op},k}(B_E,1_BX)$ denotes the $\Tor$ groups of the
$E\ot B^{\op}$-module $B_E$, relative to the family of the $E\ot
B^{\op}$-linear epimorphisms which split as $k$-linear morphisms.
\endproclaim

\remark{Remark} When the present paper was finished we learned that
Theorem~1 was also obtained in \cite{C-M-R-S} under the additional
assumptions that $k$ is a field, $X=E$ and $M$ is $A$-projective on the
left or $B$-projective on the right.
\endremark

\head Proof of the results\endhead

Let $(E^{*+2},b'_*)$ be the canonical resolution of $E$ and let
$(X_*,b'_*)$ be the $E$-bimodule subcomplex of  $(E^{*+2},b'_*)$, defined
by
$$
X_n = A^{n+2} \oplus B^{n+2} \oplus \bigoplus_{i=0}^{n+1} A^i \ot M \ot
B^{n+1-i}.
$$
It is easy to see that $(X_*,b'_*)$ is a direct summand of $(E^{*+2},b'_*)$
as an $E$-bimodule complex. Moreover, the complex
$$
E @<b'_0<< X_0 @<b'_1<< X_1 @<b'_2<< X_2 @<b'_3<< X_3 @<b'_4<< X_4 @<b'_5<<
X_5 @<b'_6<< X_6 @<b'_7<< X_7 @<b'_8<< \dots
$$
is contractible as a right $E$-module complex. Hence, $(X_*,b'_*)$ is a
projective resolution of the $E^e$-module $E$, relative to the family of
the $E^e$-linear epimorphisms which split as $k$-linear morphisms.

\smallskip

Let $(X_*^A,b'_*)$ and $(X_*^B,b'_*)$ be the subcomplexes of $(X_*,b'_*)$,
defined by $X_n^A = A^{n+1}\ot (A\oplus M)$ and $X_n^B = (B\oplus M) \ot
B^{n+1}$. It is easy to see that $(X_*^A,b'_*)$ and $(X_*^B,b'_*)$ are
projective resolutions of the $E^e$-modules $1_AE$ and $E1_B$ respectively,
relative to the family of the $E^e$-linear epimorphisms which split as
$k$-linear morphisms. We have the following:

\proclaim{Lemma 3} Let $\mu\: \tfrac{X_1}{X_1^A\oplus X_1^B}\to M$ be the
map defined by $\mu(a\ot m \ot b) = amb$, for $a\in A$, $b\in B$ and $m\in
M$. The complex
$$
M @<\mu<< \frac{X_1}{X_1^A\oplus X_1^B} @<b'_2<< \frac{X_2}{X_2^A\oplus
X_2^B} @<b'_3<< \frac{X_3}{X_3^A\oplus X_3^B} @<b'_4<<
\frac{X_4}{X_4^A\oplus X_4^B} @<b'_5<< \dots,\tag*
$$
is a relative projective resolution of $M$ as an $E$-bimodule. A
contracting homotopy of $(*)$ as a complex of $k$-modules is the the family
$\si_1\:M \to \tfrac{X_1}{X_1^A \oplus X_1^B}$ and
$\si_{n+1}\:\tfrac{X_n}{X_n^A\oplus X_n^B} \to
\tfrac{X_{n+1}}{X_{n+1}^A\oplus X_{n+1}^B}$ ($n\ge 1$), defined by:
$$
\align
&\si_1(m) = 1_A\ot m \ot 1_B,\\
&\si_{n+1}(a_0\ot m\ot \bb_{2,n+1}) = 1_A\ot a_0\ot m\ot \bb_{2,n+1} +
(-1)^n 1_A\ot a_0m \ot \bb_{2,n+1} \ot 1_B,\\
&\si_{n+1}(\ba_{0,i}\ot m\ot \bb_{i+2,n+1}) = 1_A\ot \ba_{0,i}\ot m\ot
\bb_{i+2,n+1} \quad \text{for $i>0$,}
\endalign
$$
where $\ba_{0i} = a_0\ot\cdots\ot a_i$ and $\bb_{i+2,n+1} =
b_{i+2}\ot\cdots\ot b_{n+1}$.
\endproclaim

\demo{Proof} It follows by a direct computation.\qed
\enddemo

\proclaim{Lemma 4} We have
$$
\align
& \Hom_{A^e}((A^{*+2},b'_*),X_{AA}) \simeq \Hom_{E^e}((X^A_*,b'_*),X),\\
& \Hom_{B^e}((B^{*+2},b'_*),X_{BB}) \simeq \Hom_{E^e}((X^B_*,b'_*),X).
\endalign
$$
\endproclaim

\demo{Proof} Since, for every $f\in \Hom_{A^e}(A^{n+2},X)$,
$$
f(a_0\ot\cdots \ot a_{n+1}) = \Bigl(\smallmatrix 1 & 0 \\ 0 &
0\endsmallmatrix \Bigr)f(a_0\ot\cdots \ot a_{n+1})\Bigl(\smallmatrix 1 & 0
\\ 0 & 0\endsmallmatrix \Bigr)\in X_{AA},
$$
the canonical inclusion $i_n\: \Hom_{A^e}(A^{n+2},X_{AA})\to
\Hom_{A^e}(A^{n+2},X)$ is an isomorphism. Let $\theta_n^A \:
\Hom_{A^e}(A^{n+2},X) \to \Hom_{E^e}(X^A_n,X)$ be the map defined by
$$
\alignat2
& \theta_n^A(f)(a_0\ot\cdots\ot a_{n+1}) = f(a_0\ot\cdots \ot a_{n+1}) &&
\quad\text{for $a_i\in A$,}\\
& \theta_n^A(f)(a_0\ot\cdots\ot a_n\ot m) = f(a_0\ot\cdots\ot a_n\ot 1_A)
\Bigl(\smallmatrix 0 & m \\ 0 & 0\endsmallmatrix \Bigr) && \quad\text{for
$a_i\in A$ and $m\in M$,}
\endalignat
$$
and let $\vartheta_n^A \: \Hom_{E^e}(X^A_n,X) \to \Hom_{A^e}(A^{n+2},X)$ be
the map defined by restriction. Clearly $\vartheta_n^A \circ \theta_n^A =
id$. Let us see that $\theta_n^A \circ \vartheta_n^A = id$. Let $\varphi\in
\Hom_{E^e}(X^A_n,X)$. It is clear that $\theta_n^A \circ \vartheta_n^A
(\varphi)(a_0\ot\cdots\ot a_{n+1}) = \varphi(a_0\ot\cdots\ot a_{n+1})$ for
all $a_0,\dots,a_{n+1}\in A$. Since
$$
\align
\varphi(a_0\ot\cdots\ot a_n\ot m) & = \varphi(a_0\ot\cdots\ot a_n\ot 1_A)
\Bigl(\smallmatrix 0 & m \\ 0 & 0\endsmallmatrix \Bigr)\\
& = \theta_n^A(\vartheta_n^A(\varphi))(a_0\ot\cdots\ot a_n\ot
1_A)\Bigl(\smallmatrix 0 & m \\ 0 & 0\endsmallmatrix \Bigr)\\
& = \theta_n^A(\vartheta_n^A(\varphi))(a_0\ot\cdots\ot a_n\ot m),
\endalign
$$
for all $a_0,\dots,a_n\in A$ and $m\in M$, we have that $\theta_n^A \circ
\vartheta_n^A (\varphi) = \varphi$. As the family $\theta_*\circ i_*$ is a
map of complexes, the first assertion holds. The proof of the second one is
similar.\qed
\enddemo

\proclaim{Lemma 5} We have
$$
\Hom_{A\ot B^{\op}}\bigl(\bigl(\tfrac{X_*} {X_*^A\oplus
X_*^B},b'_*\bigr),X_{AB} \bigr) \simeq \Hom_{E^e}\bigl(\bigl(
\tfrac{X_*}{X_*^A \oplus X_*^B},b'_*\bigr),X\bigr).
$$
\endproclaim

\demo{Proof} Since, for every $f\in \Hom_{E^e}\bigl(\tfrac{X_n}{X_n^A \oplus
X_n^B},X\bigr)$,
$$
f(x_0\ot\cdots \ot x_{n+1}) = \Bigl(\smallmatrix 1 & 0 \\ 0 &
0\endsmallmatrix \Bigr) f(x_0\ot\cdots \ot x_{n+1})\Bigl(\smallmatrix 0 & 0
\\ 0 & 1\endsmallmatrix \Bigr) \in X_{AB},
$$
the canonical inclusion $\Hom_{E^e}\bigl(\tfrac{X_n}{X_n^A\oplus
X_n^B},X_{AB}\bigr)\to \Hom_{E^e}\bigl( \tfrac{X_n} {X_n^A \oplus
X_n^B},X\bigr)$ is an isomorphism. To end the proof it suffices to observe
that
$$
\Hom_{E^e}\bigl(\tfrac{X_n}{X_n^A\oplus X_n^B},X_{AB}\bigr)\simeq
\Hom_{A\ot B^{\op}}\bigl(\tfrac{X_n}{X_n^A\oplus X_n^B},X_{AB}\bigr).\qed
$$
\enddemo

\subheading{Proof of Theorem 1} Because of Lemma~3, the short exact sequence
$$
0 @>>> (X_*^A,b'_*)\oplus (X_*^B,b'_*) @>>> (X_*,b'_*) @>>>
\left(\frac{X_*}{X_*^A\oplus X_*^B},b'_*\right) @>>> 0,
$$
gives rise to the long exact sequence
$$
\align
0 & @>>> \Ext^0_{E^e,k}(E,X) @>>> \Ext^0_{E^e,k}(1_AE\oplus E1_B,X) @>>>
\Ext^0_{E^e,k}(M,X) @>>> \\
& @>>> \Ext^1_{E^e,k}(E,X) @>>> \Ext^1_{E^e,k}(1_AE\oplus E1_B,X) @>>>
\Ext^1_{E^e,k}(M,X) @>>>\dots.
\endalign
$$
To end the proof it suffices to apply Lemmas~4 and 5.\qed

\medskip

\proclaim{Lemma 6} Let $\mu\: \tfrac{X_0}{X_0^A}\to B_E$ be the
map defined by $\mu(b_0\ot b_1+ m\ot b) = b_0b_1$, for $b,b_0,b_1\in B$ and
$m\in M$. The complex
$$
B_E @<\mu<< \frac{X_0}{X_0^A} @<b'_1<< \frac{X_1}{X_1^A} @<b'_1<<
\frac{X_2}{X_2^A} @<b'_3<< \frac{X_3}{X_3^A} @<b'_4<< \frac{X_4}{X_4^A}
@<b'_5<< \frac{X_5}{X_5^A} @<b'_6<< \dots,\tag*
$$
is a relative projective resolution of $B_E$ as an $E$-bimodule. A
contracting homotopy of $(*)$ as a complex of $k$-modules is the family
$\si_0\:B_E \to \tfrac{X_0}{X_0^A}$ and $\si_{n+1}\:\tfrac{X_n}{X_n^A} \to
\tfrac{X_{n+1}}{X_{n+1}^A}$ ($n\ge 0$), defined by:
$$
\si_{n+1}(x_0\ot \cdots \ot x_n) = 1_A\ot x_0\ot \cdots \ot x_n.
$$
\endproclaim

\demo{Proof} It follows by a direct computation.\qed
\enddemo

\proclaim{Lemma 7} We have $\Hom_{E\ot B^{\op}}\bigl(\bigl(\tfrac{X_*}
{X_*^A},b'_*\bigr),X1_B\bigr)\simeq \Hom_{E^e}\bigl(\bigl(\tfrac{X_*}
{X_*^A},b'_*\bigr),X\bigr)$.
\endproclaim

\demo{Proof} Since, for every $f\in \Hom_{E^e}\bigl(\tfrac{X_n}{X_n^A},
X\bigr)$,
$$
f(x_0\ot\cdots \ot x_{n+1}) = f(x_0\ot\cdots\ot x_{n+1})\Bigl(\smallmatrix
0 & 0 \\ 0 & 1\endsmallmatrix \Bigr) \in X1_B,
$$
the canonical inclusion $\Hom_{E^e}\bigl(\tfrac{X_n}{X_n^A},X1_B\bigr)\to
\Hom_{E^e}\bigl(\tfrac{X_n}{X_n^A},X\bigr)$ is an isomorphism. To end the
proof it suffices to observe that
$$
\Hom_{E\ot B^{\op}}\bigl(\tfrac{X_n} {X_n^A},X1_B\bigr) =
\Hom_{E^e}\bigl(\tfrac{X_n}{X_n^A},X1_B\bigr).\qed
$$
\enddemo

\subheading{Proof of Theorem 2} Because of Lemma~6, the short exact sequence
$$
0 @>>> (X_*^A,b'_*) @>>> (X_*,b'_*) @>>> \left(\frac{X_*}{X_*^A},
b'_*\right) @>>> 0,
$$
gives rise to the long exact sequence
$$
\align
0 & @>>> \Ext^0_{E^e,k}(B_E,X) @>>> \Ext^0_{E^e,k}(E,X) @>>>
\Ext^0_{E^e,k}(1_AE,X) @>>> \\
& @ >>> \Ext^1_{E^e,k}(B_E,X) @>>> \Ext^1_{E^e,k}(E,X) @>>>
\Ext^1_{E^e,k}(1_AE,X) @>>> \dots.
\endalign
$$
To end the proof it suffices to apply Lemmas~4 and 7.\qed

\medskip

\proclaim{Lemma 8} Let $((A\ot B^{\op})^{*+1}\ot M),b'_*)$ be the bar
resolution of $M$ as a left $A\ot B^{\op}$-module. There is a map of
resolutions $\ga_*\:((A\ot B^{\op})^{*+1}\ot M),b'_*) \to \left(
\frac{X_{*+1}}{X_{*+1}^A\oplus X_{*+1}^B},b'_{*+1}\right)$, defined by
$$
\ga_n((a_0\ot b_0)\ot\cdots\ot (a_n\ot b_n)\ot m)) = \sum_{i=0}^n
(-1)^{\binom{n-i}{2}} \ba_{0,i}\ot \fa_{i+1,n}m\ot \bb_{n,i+1}\ot
\fb_{i,0},
$$
where
$$
\alignat2
&\ba_{0,i} = a_0\ot a_1\ot\cdots \ot a_i,&&\qquad \fa_{i+1,n} =
a_{i+1}a_{i+2}\cdots a_{n},\\
&\bb_{n,i+1} = b_n\ot b_{n-1}\ot\cdots \ot b_{i+1},&&\qquad \fb_{i,0} =
b_ib_{i-1}\cdots b_0.
\endalignat
$$
\endproclaim

\demo{Proof} We must prove that $\mu \circ \ga_0 = b'_0$, where $\mu$ is
the map introduced in Lemma~3 and that $\ga_{n-1}\circ b'_n = b'_n \circ
\ga_n$ for all $n\ge 1$. The first assertion is evident. Let us check the
second one. We have
$$
\allowdisplaybreaks
\align
& \ga_{n-1}\circ b'_n ((a_0\ot b_0)\ot\cdots\ot (a_n\ot b_n)\ot m)\\
& = \ga_{n-1}\biggl(\sum_{j=0}^{n-1} (-1)^j (a_0\ot b_0)\ot\cdots\ot
(a_ja_{j+1}\ot b_{j+1}b_j)\ot\cdots \ot (a_n\ot b_n)\ot m\biggr)\\
& + (-1)^n \ga_{n-1}((a_0\ot b_0)\ot\cdots \ot (a_{n-1}\ot b_{n-1})\ot
a_nmb_n) \\
& = \sum_{j=0}^{n-1} \sum_{i=j}^{n-1} (-1)^{j+\binom{n-i-1}{2}} \ba_{0,j-1}
\ot a_ja_{j+1}\ot\ba_{j+2,i+1}\ot\fa_{i+2,n}m\ot\bb_{n,i+2}\ot\fb_{i+1,0}\\
& + \sum_{j=0}^{n-1} \sum_{i=0}^{j-1} (-1)^{j+\binom{n-i-1}{2}} \ba_{0,i}\ot
\fa_{i+1,n}m\ot \bb_{n,j+2}\ot b_{j+1}b_j\ot\bb_{j-1,i-1}\ot \fb_{i,0} \\
& + \sum_{i=0}^{n-1} (-1)^{n+\binom{n-i-1}{2}} \ba_{0,i}\ot \fa_{i+1,n}mb_n
\ot \bb_{n-1,i+1}\ot \fb_{i,0}\\
& = \sum_{i=1}^n \sum_{j=0}^{i-1} (-1)^{\binom{n-i}{2}+j} \ba_{0,j-1}\ot
a_ja_{j+1}\ot \ba_{j+2,i}\ot \fa_{i+1,n}m\ot \bb_{n,i+1}\ot \fb_{i,0}\\
& + \sum_{i=0}^{n-2} \sum_{j=i+1}^{n-1} (-1)^{\binom{n-i}{2}+i+1+n-j} \ba_{0,i}
\ot\fa_{i+1,n}m\ot \bb_{n,j+2}\ot b_{j+1}b_j\ot\bb_{j-1,i}\ot \fb_{i-1,0}\\
& + \sum_{i=0}^{n-1} (-1)^{\binom{n-i}{2}+i+1}\ba_{0,i}\ot \fa_{i+1,n}mb_n
\ot \bb_{n-1,i+1}\ot \fb_{i,0}\\
& = b'_n\left(\sum_{i=0}^n (-1)^{\binom{n-i}{2}} \ba_{0,i}\ot
\fa_{i+1,n}m\ot \bb_{n,i+1}\ot \fb_{i,0} \right)\\
& = b'_n \circ \ga_n ((a_0\ot b_0)\ot\cdots\ot (a_n\ot b_n)\ot m).\qed
\endalign
$$
\enddemo

\subheading{Proof of Theorem 3} 1) It is easy to see that $\H^*(E,E)
@>j^*>> \H^*(A,A) \oplus \H^*(B,B)$ is induced by the canonical restriction
$$
\Hom_{E^e}((E^{*+2},b'_*),E) @>>>\Hom_{A^e}((A^{*+2},b'_*),A) \oplus
\Hom_{B^e}((B^{*+2},b'_*),B).
$$
From this fact follows immediately that $j^*$ is a map of graded rings.

\smallskip

\noindent 2) We prove the second assertion. The first one follows
similarly. It is easy to see that $\H^*(B,B) @>>> \Ext^*_{A\ot
B^{\op},k}(M,M)$ is induced by the map of complexes
$$
\Hom_{B^e}((B^{*+2},b'_*),B) @>\wt{\de}^*>> \Hom_{A\ot B^{\op}}\bigl(\bigl(
\tfrac{X_{*+1}}{X_{*+1}^A\oplus X_{*+1}^B},-b'_{*+1}\bigr),M\bigr),
$$
defined by
$$
\wt{\de}^n(f)(\ba_{0,i}\ot m\ot \bb_{i+1,n+1}) = \cases (-1)^n
f(\ba_{0,n}\ot 1_A)mb_{n+1}& \text{ if $i=n$,}\\ 0 & \text{ in other
cases,} \endcases
$$
where $f\in \Hom_{B^e}(B^{n+2},B)$, $\ba_{0,i} = a_0\ot\cdots\ot a_i$ and
$\bb_{i+1,n+1} = b_{i+1}\ot\cdots\ot b_{n+1}$. Let us consider the morphism
$$
\Hom_{A\ot B^{\op}}\bigl(\bigl(\tfrac{X_{*+1}}{X_{*+1}^A\oplus
X_{*+1}^B},b'_{*+1}\bigr),M\bigr) @>\wt{\ga}^*>> \Hom_{A\ot
B^{\op}}\bigl(((A\ot B^{\op})^{*+1}\ot M),b'_*),M\bigr)
$$
induced by the map $\ga_*$ of Lemma~8. Let $\phi^* = \wt{\ga}^*\circ
\wt{\de}^*$. It is immediate that
$$
\phi^n(f)((a_0\ot b_0)\ot\cdots\ot (a_n\ot b_n)\ot m) = (-1)^{\binom{n}{2}}
a_0\cdots a_n m f(1_B\ot b_n\ot \cdots\ot b_0).
$$
Hence, for $f\in \Hom_{B^e}(B^{r+2},B)$, $g\in \Hom_{B^e}(B^{s+2},B)$
and $x_0\ot\cdots \ot x_{r+s} = (a_0\ot b_0)\ot\cdots\ot (a_{r+s}\ot
b_{r+s})$, we have
$$
\align
&\phi^r(f)\phi^s(g)(x_0\ot\cdots \ot x_{r+s}\ot m)\\
& = \phi^r(f) (x_0\ot\cdots \ot x_r\ot \phi^s(g)((1_A\ot 1_B)\ot
x_{r+1}\ot\cdots \ot x_{r+s}\ot m)\\
& = (-1)^{\binom{s}{2}}\phi^r(f) (x_0\ot\cdots \ot x_r\ot a_{r+1}\cdots
a_{r+s} m g(1_B\ot b_{r+s}\ot\cdots \ot b_{r+1}\ot 1_B))\\
& = (-1)^{rs+\binom{r+s}{2}} a_0\cdots a_{r+s}m g(1_B\ot
b_{r+s}\ot\cdots \ot b_{r+1}\ot 1_B) f(1_B\ot b_r\ot\cdots \ot b_0)\\
& = (-1)^{rs}(-1)^{\binom{r+s}{2}} a_0\cdots a_{r+s}m (gf)(1_B\ot
b_{r+s}\ot\cdots \ot b_0)\\
& = \phi^{r+s}((-1)^{rs}gf)(x_0\ot\cdots \ot x_{r+s}\ot m).
\endalign
$$
This finished the proof, since $\H^*(B,B)$ is graded commutative.

\smallskip

\noindent 3) The complex $\Hom_{E^e}((X_*,b'_*),E)$ is a differential
graded algebra with the product defined by $(fg)(x_0\ot\cdots\ot x_{r+s+1})
= f(x_0\ot\cdots\ot x_r\ot 1_E)g(1_E\ot x_{r+1}\ot\cdots\ot x_{r+s+1})$,
for $f\in \Hom_{E^e}(X_r,E)$ and $g\in \Hom_{E^e}(X_s,E)$. It is immediate
that this product induces the cup product in $\H^*(E,E)$. Assume that $s\ge
1$ and that $f$ belongs to the image of $\Hom_{E^e}\left(\frac{X_r}{X_r^A
\oplus X_r^B},E \right) @>>> \Hom_{E^e}(X_r,E)$. Let $g'\in
\Hom_{E^e}(X_s,E)$ the cocycle defined by
$$
g'(x_0\ot\cdots\ot x_r) = \cases g(x_0\ot\cdots\ot x_r) &\text{if
$x_0\ot\cdots\ot x_r\in A^r\ot (A\oplus M)$} \\ 0 &\text{in other case}
\endcases
$$
It is easy to check that $gf = g'f$ and that $fg' = 0$. Since $g'f$ is
homologous to $(-1)^{rs}fg'$ we have that the class of $gf$ in
$\H^{r+s}(E,E)$ is zero.\qed

\enddocument